\documentclass[11pt, reqno]{amsart}

\usepackage[OT2,T1]{fontenc}

\usepackage[usenames,dvipsnames]{color}

\usepackage{url}
\usepackage{amsmath}
\usepackage{array}
\usepackage{graphicx}
\usepackage{amssymb}
\usepackage{amsthm}
\definecolor{link}{RGB}{11,0,128}
\usepackage[colorlinks=true,citecolor=NavyBlue,linkcolor=Bittersweet,urlcolor=link]{hyperref}
\usepackage{paralist}
\usepackage{colonequals}		
\usepackage{color}
\usepackage{mathabx}		
\usepackage{amscd}
\usepackage[all,cmtip]{xy}	
\usepackage{sseq}			
\usepackage{verbatim}
\usepackage{parskip}
\usepackage{microtype}
\usepackage[in]{fullpage}
\usepackage{mathrsfs} 			
\usepackage{cleveref}
\usepackage{enumitem}
\usepackage{moreenum}			
\usepackage[alphabetic,lite]{amsrefs} 	
\usepackage{stmaryrd}	
\usepackage[foot]{amsaddr}
\usepackage{ragged2e}

\DeclareSymbolFont{cyrletters}{OT2}{wncyr}{m}{n}
\DeclareMathSymbol{\Sha}{\mathalpha}{cyrletters}{"58}

\newcommand{\gA}{\alpha}
\newcommand{\gB}{\beta}

\newcommand{\bF}{\mathbb{F}}

\newcommand{\bP}{\mathbb{P}}
\newcommand{\bQ}{\mathbb{Q}}

\newcommand{\bZ}{\mathbb{Z}}

\newcommand{\bbB}{\mathbf{B}}

\newcommand{\cE}{\mathcal{E}}

\newcommand{\cG}{\mathcal{G}}

\newcommand{\sU}{\mathscr{U}}
\newcommand{\sV}{\mathscr{V}}

\newcommand{\sX}{\mathscr{X}}
\newcommand{\sY}{\mathscr{Y}}

\newcommand{\ra}{\rightarrow}

\newcommand{\wt}{\widetilde}
\newcommand{\wh}{\widehat}

\newcommand{\ce}{\colonequals}

\renewcommand{\b}{\textbf}
\newcommand{\surjects}{\twoheadrightarrow}

\newcommand{\st}{{\mathrm{st}}}		

\providecommand{\p}[1]{\left(#1\right)}

\providecommand{\f}[2]{\frac{#1}{#2}}
\DeclareMathOperator{\Ker}{Ker}			
\DeclareMathOperator{\Ver}{Ver}	
\DeclareMathOperator{\GL}{GL}		
\DeclareMathOperator{\Frob}{Frob}		
\newcommand{\ba}{\begin{aligned}}
\newcommand{\ea}{\end{aligned}}
\newcommand{\be}{\begin{equation}}
\newcommand{\ee}{\end{equation}}
\newcommand{\pf}{\begin{proof}}
\newcommand{\bpf}{\begin{proof}}
\newcommand{\epf}{\end{proof}}
\newcommand{\bthm}{\begin{thm}}
\newcommand{\ethm}{\end{thm}}

\newcommand{\bprop}{\begin{prop}}
\newcommand{\eprop}{\end{prop}}
\newcommand{\bcor}{\begin{cor}}
\newcommand{\ecor}{\end{cor}}

\newcommand{\brem}{\begin{rem}}
\newcommand{\erem}{\end{rem}}

\newcommand{\brems}{\begin{rems} \hfill \begin{enumerate}[label=\b{\thesubsection.},ref=\thesubsection]}
\newcommand{\bremstweak}{\begin{rems-tweak} \hfill \begin{enumerate}[label=\b{\thesubsection.},ref=\thesubsection]}
\newcommand{\bremst}{\begin{rems-tweak} \hfill \begin{enumerate}[label=\b{\thesubsection.},ref=\thesubsection]}
\newcommand{\remi}{\addtocounter{subsection}{1} \item}

\newcommand{\erems}{\end{enumerate} \end{rems}}
\newcommand{\eremstweak}{\end{enumerate} \end{rems-tweak}}
\newcommand{\eremst}{\end{enumerate} \end{rems-tweak}}
\newcommand{\blem}{\begin{lemma}}
\newcommand{\elem}{\end{lemma}}
\newcommand{\blemt}{\begin{lemma-tweak}}
\newcommand{\elemt}{\end{lemma-tweak}}
\newcommand{\bconj}{\begin{conj}}
\newcommand{\econj}{\end{conj}}
\newcommand{\bprob}{\begin{Problem}}
\newcommand{\eprob}{\end{Problem}}

\newcommand{\bq}{\begin{q}}
\newcommand{\eq}{\end{q}}
\newcommand{\benum}{\begin{enumerate}[label={(\alph*)}]}
\newcommand{\benuma}{\begin{enumerate}[label={(\arabic*)}]}
\newcommand{\benumr}{\begin{enumerate}[label={(\roman*)}]}
\newcommand{\eenum}{\end{enumerate}}
\newcommand{\bc}{}
\newcommand{\bd}{\begin{defn}}
\newcommand{\ed}{\end{defn}}

\newcommand{\beg}{\begin{eg}}
\newcommand{\eeg}{\end{eg}}

\newcommand{\bcl}{\begin{claim}}
\newcommand{\ecl}{\end{claim}}
\newcommand{\lab}{\label}

\newcommand{\q}{\quad}
\newcommand{\qq}{\quad\quad}

\newcommand{\tst}{\textstyle}

\newcommand{\bal}{\mathrm{bal}}

\newcommand*{\QED}{\hfill\ensuremath{\qed}}

\theoremstyle{plain}
\newtheorem{thm}[subsection]{Theorem}
\Crefname{thm}{Theorem}{Theorems}

\newtheorem{thm-tweak}[subsection]{Theorem}
\Crefname{thm-tweak}{Theorem}{Theorems}

\Crefname{rethm}{Theorem}{Theorem}
\newtheorem{prop}[subsection]{Proposition}
\Crefname{prop}{Proposition}{Propositions} 

\newtheorem{prop-tweak}[subsection]{Proposition}
\Crefname{prop-tweak}{Proposition}{Propositions} 

\newtheorem{Q}[subsection]{Question}
\Crefname{Q}{Question}{Questions}
\Crefname{eg}{Example}{Examples}
\newtheorem{Problem}[subsection]{Problem}
\Crefname{Problem}{Problem}{Problems}
\newtheorem{conj}[subsection]{Conjecture}
\Crefname{conj}{Conjecture}{Conjectures}
\newtheorem{cor}[subsection]{Corollary}
\Crefname{cor}{Corollary}{Corollaries}
\newtheorem{cor-tweak}[subsection]{Corollary}
\Crefname{cor-tweak}{Corollary}{Corollaries}

\newtheorem{lemma}[subsection]{Lemma}

\newtheorem{lemma-tweak}[subsection]{Lemma}
\Crefname{lemma-tweak}{Lemma}{Lemmas}

\Crefname{subprop}{Proposition}{Propositions}
\Crefname{subcor}{Corollary}{Corollaries}
\Crefname{sublem}{Lemma}{Lemmas}

\theoremstyle{remark}
\newtheorem{claim}[equation]{Claim}
\Crefname{claim}{Claim}{Claims}

\Crefname{subrem}{Remark}{Remarks}

\theoremstyle{definition}
\newtheorem{defn}[subsection]{Definition}
\Crefname{defn}{Definition}{Definitions}
\newtheorem{defn-tweak}[subsection]{Definition}
\Crefname{defn-tweak}{Definition}{Definitions}

\Crefname{conv}{Convention}{Conventions}
\newtheorem{eg}[subsection]{Example}
\newtheorem{rem}[subsection]{Remark}
\Crefname{rem}{Remark}{Remarks}

\newtheorem{eg-tweak}[subsection]{Example}
\Crefname{eg-tweak}{Example}{Examples}

\newtheorem{rem-tweak}[subsection]{Remark}
\Crefname{rem-tweak}{Remark}{Remarks}

\newtheorem*{rems}{Remarks}

\newtheorem*{rems-tweak}{Remarks}

\newtheoremstyle{subsection-tweak}
   {11pt}
   {3pt}%
   {}
   {}%
   {\bfseries}
   {}%
   {.5em}
   {\thmnumber{\@{#1}{}\@{#2}.}%
    \thmnote{~{\bfseries#3.}}}

\Crefname{innercustomconj}{Conjecture}{Conjecture}

\theoremstyle{subsection-tweak}
\newtheorem{pp}[subsection]{}
\newcommand{\bpp}{\begin{pp}}
\newcommand{\epp}{\end{pp}}

\theoremstyle{subsection-tweak}
\newtheorem{pp-tweak}[subsection]{}

\numberwithin{equation}{subsection}

\makeatletter
\def\@tocline#1#2#3#4#5#6#7{
    \begingroup 
    \@ifempty{#4}{%
    }{%
    }%

    \parindent\z@ \leftskip#3\relax \advance\leftskip\@tempdima\relax
    #5\hskip-\@tempdima
      \ifcase #1
       \or\or \hskip 2em \or \hskip 1em \else \hskip 3em \fi%
      #6\nobreak\relax
    \dotfill\hbox to\@pnumwidth{\@tocpagenum{#7}}\par
    \nobreak
    \endgroup
  }
 \def\l@section{\@tocline{1}{0pt}{1pc}{}{}}

\renewcommand{\tocsection}[3]{%
  \indentlabel{\@ifnotempty{#2}{\makebox[1.3em][l]{%
    \ignorespaces#1 \bfseries{#2}.\hfill}}}\bfseries{#3}
    \vspace{1.5pt}}

\renewcommand{\tocsubsection}[3]{%
  \indentlabel{\@ifnotempty{#2}{\hspace*{-0.5em}\makebox[2.1em][l]{%
    \ignorespaces#1#2.\hfill}}}#3
    \vspace{1.5pt}}

\makeatother

\begin{document}
\author{K\k{e}stutis \v{C}esnavi\v{c}ius}
\title{Coarse base change fails for some modular curves}
\date{\today}
\address{Department of Mathematics, University of California, Berkeley, CA 94720-3840, USA}
\email{kestutis@berkeley.edu}

\begin{abstract} 
For a congruence level $H \subset \GL_2(\wh{\bZ})$, the formation of the modular curve $X_H$, i.e.,~of the coarse moduli space of the level $H$ modular stack $\sX_H$, is known to commute with arbitrary base change in a wide range of cases. We exhibit infinitely many $H$, for instance, $H = \Gamma_1(4)$, for which this coarse base change property fails. In our examples failure is witnessed for base change to $\bF_2$ and for any $\bZ_{(2)}$-fiberwise dense open substack of $\sX_H$. These examples fill in several open entries in a table in the book of Katz and Mazur.
\end{abstract}

\maketitle

\section{Introduction}

\bpp[Coarse base change for modular curves]
In the study of modular forms modulo a prime $p$, one often considers the $\bF_p$-fiber $(X_H)_{\bF_p}$ of the modular curve $X_H$ over $\bZ$ for an open subgroup $H \subset \GL_2(\wh{\bZ})$. This $\bF_p$-fiber may be difficult to analyze---$X_H$ is defined as the coarse moduli space of the level $H$ modular stack $\sX_H$, but there is no a priori reason why $(X_H)_{\bF_p}$ should be the coarse space of $(\sX_H)_{\bF_p}$ if $p = 2$ or $p = 3$. If, however, $(X_H)_{\bF_p}$ happens to be the coarse space of $(\sX_H)_{\bF_p}$ for every prime $p$, then one knows (say, from \cite{Ces15a}*{3.3.1}) that \emph{coarse base change} holds for $\sX_H$, i.e.,~that $(X_H)_S$ is the coarse space of $(\sX_H)_S$ for every scheme $S$. For instance, this is the case when $H = \GL_2(\wh{\bZ})$.  The main goal of this paper is to show that coarse base change fails for some~$\sX_H$.
\epp

\bpp[Known results on coarse base change for modular curves] \lab{known}
To put our examples in perspective, we turn to known results, all of which confirm coarse base change for $\sX_H$ for suitable~$H$.

Only primes $2$ and $3$ may divide the order of the automorphism group of an elliptic curve over a field, so it is a generality that $(X_H)_S$ is the coarse moduli space of $(\sX_H)_S$ whenever $S$ is a $\bZ[\f{1}{6}]$-scheme (the essential observation is that the formation of the ring of invariants $R^G$ for the action of a finite group $G$ on a ring $R$ commutes with arbitrary base change whenever $\#G \in R^\times$). As is explained in \cite{Ces15a}*{6.4~(b)}, it then follows from the results of Deligne--Rapoport and Katz--Mazur that the same holds whenever $S$ is a $\bZ[\f{1}{\gcd(6, n)}]$-scheme for any $n \in \bZ_{\ge 1}$ with $\Gamma(n) \subset H$, where $\Gamma(n) \ce \Ker(\GL_2(\wh{\bZ}) \surjects \GL_2(\bZ/n\bZ))$. 

There are also situations in which $2$ or $3$ divides the level and coarse base change continues to hold (with $\sX_H \not\cong X_H$). For instance, by \cite{DR73}*{VI.6.10} (with \cite{Ces15a}*{3.3.1} as before), this is the case when $H = \Gamma_0(6) \cap H'$ with an $H'$ for which $\Gamma(n) \subset H'$ for some $n \in \bZ_{\ge 1}$ that is prime to $6$.

One may also ask whether coarse base change holds generically, say, whether it holds for the open substack $\sU \subset \sX_H$ on which the $j$-invariant avoids the set $\{0, 1728, \infty\}$ (this automatically disposes of the need to consider the $\bF_3$-fiber). To address this, Katz and Mazur compiled a table \cite{KM85}*{8.5.4} that confirms coarse base change for $\sU$ for many common choices of $H$. The examples of this paper fill in a ``NO'' for some open entries of this table.
\epp

\bpp[The nature of our examples]
We choose the subgroup $H$ to be either
\be \lab{which-H}
H = \Gamma_1(4) \cap \Gamma_0(4n) \q \text{for an $n \in \bZ_{\ge 1}$} \qq \text{or} \qq H = \Gamma_1^\bal(4)
\ee
(see \S\S\ref{X1(4)}--\ref{bal} for a review of the definitions of these subgroups) and show that $(X_H)_{\bF_2}$ is not the coarse moduli space of $(\sX_H)_{\bF_2}$. For the sake of clarity, we first treat the simplest case $H = \Gamma_1(4)$ in \S\ref{X1(4)} and then refine the argument in \S\ref{G1-G0} to include all $H$ of the form $\Gamma_1(4) \cap \Gamma_0(4n)$. This family of examples shows that coarse base change may fail for modular curves of arbitrarily high genus.

In our examples the failure of coarse base change persists to fiberwise dense open substacks (such as $\sU$ of \S\ref{known}) and results from the jumping of generic stabilizers: for $H$ as in \eqref{which-H}, $(\sX_H)_\bQ$ is a scheme whereas some open substack of $(\sX_H)_{\bF_2}$ is a $\bZ/2\bZ$-gerbe over its coarse space. The conclusion then results from the fact that over a fiberwise dense open of $(X_H)_{\bZ_{(2)}}$ the coarse moduli space morphism 
\[
(\sX_H)_{\bZ_{(2)}} \ra (X_H)_{\bZ_{(2)}}
\] 
is flat (and hence, morally, finite locally free of rank $1$), as is ensured by the miracle flatness theorem combined with an openness of the regular locus result of Nagata (in the form of \cite{EGAIV2}*{6.12.6})---see \S\ref{main}, esp.~the proof of \Cref{main-crit}, for the details of this argument.

In the view of the results of this paper, the remaining open entries in the table \cite{KM85}*{8.5.4} are those of $[\Gamma_1(N)]$ with $N = 8, 16, 32, \ldots$ and those of $[\Gamma_0(2^n; a, b)]$ with $a \ge 3$ and $b = 0$ (the assumption $a \ge b$ loses no generality and $[\Gamma_0(2^n; a, 1)] = [\Gamma_0(2^n; a, 0)]$).
\epp

It would be interesting to understand whether coarse base change may also fail at individual closed points. More precisely, it would be interesting to answer the following question.

\begin{Q}
Is there an open subgroup $H \subset \GL_2(\wh{\bZ})$ and a prime $p \in \{2, 3\}$ such that coarse base change holds for some $\bZ_{(p)}$-fiberwise dense open substack $\sV \subset (\sX_H)_{\bZ_{(p)}}$ but fails for $(\sX_H)_{\bZ_{(p)}}$ itself? In particular, can coarse base change fail for $(\sX_H)_{\bZ_{(3)}}$?
\end{Q}

\bpp[Notation and conventions]
For an open subgroup $H \subset \GL_2(\wh{\bZ})$, we let $\sX_H$ denote the level $H$ Deligne--Mumford modular stack over $\bZ$ defined in \cite{DR73}*{IV.3.3} by normalization (we also rely on the agreement discussed in \cite{Ces15a}*{6.3} with the definitions used in \cite{KM85}). Replacement of $\sX$ by $\sY$---for instance, of $\sX_H$ by $\sY_H$---indicates the elliptic curve locus. Level structures are in the sense of Drinfeld---for instance, a finite locally free $S$-subgroup $G \subset E$ of an elliptic curve $E \ra S$ is \emph{cyclic of order $n$} if fppf locally on $S$ it admits a Drinfeld $\bZ/n\bZ$-structure (the needed definitions are reviewed in \cite{Ces15a}*{4.2.2, 4.2.6, 4.2.8}). For an $n \in \bZ_{\ge 1}$, we let $\Gamma(n) \subset \GL_2(\wh{\bZ})$ be as in \S\ref{known} and we write $\sX(1)$ for $\sX_{\Gamma(1)}$ (i.e.,~for $\sX_{\GL_2(\wh{\bZ})}$). For a prime $p$, we let $(-)_{(p)}$ denote localization at $p$.
\epp

\subsection*{Acknowledgements}
I thank Brian Conrad and Barry Mazur for helpful conversations and correspondence. I thank the referee for helpful comments and suggestions. I thank the Miller Institute for Basic Research in Science at the University of California Berkeley for its support.


\section{The effect on coarse base change of the jumping of the generic stabilizer} \lab{main}

The following criterion is the source of our examples of the failure of coarse base change.

\bthm \lab{main-crit}
Suppose that $H \subset \GL_2(\wh{\bZ})$ is an open subgroup for which there exists an open substack 
\[
\sU \subset (\sY_H)_{\bZ_{(2)}}
\]
such that
\benumr
\item \lab{MC-i}
The automorphism groups of the geometric points of $\sU_\bQ$ are trivial (i.e.,~$\sU_{\bQ}$ is a scheme, cf.~\cite{Ces15a}*{4.1.4}),

\item \lab{MC-ii}
There is a nonempty open substack $\sU' \subset \sU_{\bF_2}$ whose inertia stack is the constant $\underline{\{ \pm 1\}}_{\sU'}$.
\eenum
Then the base change 
\[
\pi_{\bF_2} \colon (\sX_{H})_{\bF_2} \ra (X_H)_{\bF_2}
\] 
of the coarse moduli space morphism $\pi \colon \sX_H \ra X_H$ is not the coarse moduli space morphism of $(\sX_H)_{\bF_2}$; in fact, for the coarse moduli space $U \ce \pi_{\bZ_{(2)}}(\sU)$ of $\sU$, the map 
\[
(\pi_{\bF_2})|_{\sU_{\bF_2}}\colon \sU_{\bF_2} \ra U_{\bF_2}
\]
is not a coarse moduli space morphism.
\ethm

\brem
In practice we will have $\sU = (\sY_H)_{\bZ_{(2)}}$. However, \Cref{main-crit} with an arbitrary $\sU$ includes the claim that in presence of \ref{MC-i}--\ref{MC-ii} the failure of coarse base change for $\sX_H$ cannot be remedied by removing a $\bZ$-quasi-finite closed substack of $\sX_H$ (such as the preimage of the sections $j = 0$, $j = 1728$, and $j = \infty$ of the $j$-line $\bP^1_\bZ$).
\erem

\bpf
We replace $\sU$ by $\sU \setminus (\sU_{\bF_2} \setminus \sU')$ to be able to assume that $\sU' = \sU_{\bF_2}$. Also, it suffices to establish the last claim, so we may work over $U$ and, in particular, within the elliptic curve locus.

The flat relative $\bZ_{(2)}$-curve $U$ inherits normality from $(X_H)_{\bZ_{(2)}}$, so, by \cite{EGAIV2}*{6.12.6}, $U$ is regular away from finitely many closed points of residue characteristic $2$. We remove these points from $U$ and remove their preimages from $\sU$ to assume for the rest of the proof that $U$ is regular. 

Let $n \in \bZ_{\ge 1}$ be such that $\Gamma(n) \subset H$, let $p$ be an odd prime that does not divide $n$, and consider $\sY_{H \cap \Gamma(p)}$. The $\bZ_{(2)}$-base change $(\sY_{H \cap \Gamma(p)})_{\bZ_{(2)}}$ is a scheme because the same holds already for $(\sY_{\Gamma(p)})_{\bZ_{(2)}}$ (cf.~\cite{DR73}*{IV.2.9}). By \cite{DR73}*{proof of IV.3.9},\footnote{Even though the statement of \cite{DR73}*{IV.3.9} is incorrect (see \cite{Ces15a}*{4.5.3} for a counterexample), its proof justifies the analogous claim over the elliptic curve locus.} 
\[
\sY_{H\cap \Gamma(p)} = \sY_{H} \times_{\sY(1)} \sY_{\Gamma(p)},
\]
so, by \cite{DR73}*{the paragraph after IV.2.4},
\be \lab{deg}
(\sY_{H\cap \Gamma(p)})_{\bZ_{(2)}} \ra (\sY_{H})_{\bZ_{(2)}} \qq \text{is a $\GL_2(\bZ/p\bZ)$-torsor,}
\ee 
and hence is finite \'{e}tale of degree $\#\GL_2(\bZ/p\bZ)$. Moreover, if 
\[
\sV \subset (\sY_{H \cap \Gamma(p)})_{\bZ_{(2)}} \qq \text{denotes the preimage of}  \qq \sU \subset (\sY_{H})_{\bZ_{(2)}}
\]
(equivalently, of $U$), then the regularity of $U$ and the miracle flatness theorem \cite{EGAIV2}*{6.1.5} ensure that the finite morphism $\sV \ra U$ is locally free. Due to the $\bZ_{(2)}$-flatness of $U$, the degree of this morphism may be read off over $\bQ$, and hence, thanks to \ref{MC-i} and \eqref{deg}, equals $\#\GL_2(\bZ/p\bZ)$. For the sought conclusion, we will show that, in contrast, the map 
\[
\sV_{\bF_2} \ra U'
\]
towards the coarse moduli space $U'$ of $\sU_{\bF_2}$ is finite locally free of degree $\f{1}{2} \cdot \#\GL_2(\bZ/p\bZ)$.

Due to \ref{MC-ii}, $U'$ identifies with the rigidification $\sU' \!\!\!\!\fatslash \{\pm 1\}$ (cf.~\cite{AOV08}*{A.1} or \cite{Rom05}*{5.1}, as well as \cite{LMB00}*{8.1.1}), so $\sV_{\bF_2} \ra U'$ is finite locally free. To verify that twice its degree equals the claimed $\#\GL_2(\bZ/p\bZ)$, it remains to work \'{e}tale locally on $U'$ to first obtain an isomorphism 
\[
\sU' \simeq \bbB(\underline{\bZ/2\bZ}_{U'})
\] 
with a classifying stack and to then base change the resulting $\GL_2(\bZ/p\bZ)$-torsor 
\[
\sV_{\bF_2} \ra \bbB(\underline{\bZ/2\bZ}_{U'})
\]
along the \'{e}tale atlas $U' \ra \bbB(\underline{\bZ/2\bZ}_{U'})$ that corresponds to the trivial $\underline{\bZ/2\bZ}_{U'}$-torsor.
\epf

\brem \lab{comp-mor}
The proof shows that the comparison morphism from the coarse moduli space of $\sU'$ towards $U_{\bF_2} \cap \pi_{\bF_2}(\sU')$ is finite locally free of rank $2$ over a dense open subscheme of $U_{\bF_2} \cap \pi_{\bF_2}(\sU')$. It is a generality that this morphism is necessarily a universal homeomorphism (cf.~\cite{Ryd13}*{Thm.~6.12}).
\erem

For $H$ discussed below, the construction of suitable opens $\sU'$ rests on the following known lemmas.

\blemt[\cite{Del75}*{5.3 (III)}] \lab{Del}
If $S$ is a scheme and $E \ra S$ is an elliptic curve whose fibral $j$-invariants differ from $0$ and $1728$, then the automorphism functor of $E$ is the constant $\underline{\{ \pm 1\}}_S$. \QED
\elemt

\blemt \lab{KM}
Let $p$ be a prime and let $E \ra S$ be an elliptic curve over an $\bF_p$-scheme $S$.
\benum
\item {\upshape (\cite{KM85}*{12.2.5}).} \lab{KM-a}
For every $n \in \bZ_{\ge 1}$, the $S$-subgroup $E[p^n] \subset E$ is cyclic of order $p^{2n}$ and its standard cyclic subgroup of order $p^n$ is the kernel of the $n$-fold relative Frobenius of $E \ra S$.

\item \lab{KM-b}
If $p = 2$ and $G \subset E$ is a cyclic $S$-subgroup of order $p^2$ such that its standard cyclic subgroup $G_p$ of order $p$ is of multiplicative type and $G/G_p$ is \'{e}tale, then $G = E[p]$.
\eenum
\elemt

\bpf
The indicated reference supplies \ref{KM-a}. For \ref{KM-b}, we may use limit arguments to assume that $S$ is connected, so \cite{KM85}*{13.3.3} shows that $G$ agrees with the kernel of the composite isogeny
\[
\xymatrix@C=32pt{
E \ar[r]^-{\Frob_E} & E^{(p)} \ar[r]_{\sim}^{\iota} & \wt{E}^{(p)} \ar[r]^-{\Ver_{\wt{E}}} & \wt{E}
}
\]
for some elliptic curve $\wt{E} \ra S$ and some $S$-isomorphism $\iota$. By \Cref{Del}, any $S$-automorphism of $E^{(p)}$ is multiplication by $\pm 1$, so it remains to show that $E \simeq \wt{E}$ \'{e}tale locally on $S$, which follows from \cite{KM85}*{13.3.5 (4)} after endowing $E$ and $\wt{E}$ with compatible (\'{e}tale local on $S$) $(\bZ/3\bZ)^2$-structures (the $p = 2$ assumption ensures that the $(p -1)^\st$ infinitesimal neighborhood of the diagonal in loc.~cit.~coincides with the diagonal itself).
\epf


\section{Coarse base change fails for $\sX_1(4)$} \lab{X1(4)}

In this section we choose the level $H$ to be 
\[
\tst \Gamma_1(4) = \{ \p{ \begin{smallmatrix} a & b \\ c & d \end{smallmatrix}  } \in \GL_2(\wh{\bZ})\  |\  a \equiv 1 \bmod 4, \, c \equiv 0 \bmod 4  \}
\]
and seek to show that the formation of the coarse moduli space of $\sX_1(4) \ce \sX_{\Gamma_1(4)}$ does not commute with base change to $\bF_2$. We will obtain this from \Cref{main-crit}, so we first explain how we choose the open 
\[
\sU' \subset \sY_1(4)_{\bF_{2}}
\]
of the elliptic curve locus over $\bF_2$ (we will choose $\sU$ to be $\sY_1(4)_{\bZ_{(2)}}$).

\bpp[The choice of $\sU'$]\lab{U-prime-1}
By \cite{KM85}*{3.2} (possibly, see also \cite{Ces15a}*{4.4.4 (c)}), the Deligne--Mumford stack $\sY_1(4)$ parametrizes pairs
\[
(E \ra S,\,\gA\colon \bZ/4\bZ \ra E(S))
\]
consisting of an elliptic curve $E \ra S$ over a variable base scheme $S$ and a Drinfeld $\bZ/4\bZ$-structure $\gA$ on $E$. We consider the stack $\sY_1(4)_{\bF_2}$ together with its universal elliptic curve $\cE \ra \sY_1(4)_{\bF_2}$ and let $\cG \subset \cE$ be the finite locally free $\sY_1(4)_{\bF_2}$-subgroup of order $4$ generated by the universal $\gA$. The locus $\sU''$ of $\sY_1(4)_{\bF_2}$ over which $\cG$ is of multiplicative type is open, and we let 
\[
\sU' \subset \sU'' \qq \text{be the open locus on which the $j$-invariant satisfies} \qq j \neq 0.
\]
\epp

\bthm \lab{fail-X1(4)}
The $\bF_2$-base change of the coarse space morphism of $\sX_1(4)$ is not the coarse space morphism of $\sX_1(4)_{\bF_2}$. The same holds for any $\bZ_{(2)}$-fiberwise dense open substack of $\sX_1(4)$.
\ethm

\bpf
The claim will follow from \Cref{main-crit} once we explain why 
\[
\sU \ce \sY_1(4)_{\bZ_{(2)}} \qq \text{and} \qq \sU' \text{ of \S\ref{U-prime-1}} 
\]
satisfy its assumptions \ref{MC-i}--\ref{MC-ii}. The assumption \ref{MC-i} follows from the rigidity of $\sY_1(4)_{\bZ[\f{1}{2}]}$, i.e.,~from \cite{KM85}*{2.7.4 and 1.5.3}. For \ref{MC-ii}, we start with the nonemptiness of $\sU'$, which follows from that of $\sU''$ guaranteed by the fact that the ordinary locus of $\sY(1)_{\bF_2}$ is the forgetful image of $\sU''$ (cf.~\Cref{KM}~\ref{KM-a}).

\Cref{Del} ensures that the automorphism stack of $\cE_{\sU'}$ is the constant $\underline{\{ \pm 1\}}_{\sU'}$, so it remains to show that $[-1]_E$ preserves $\gA$ for every $(E \ra S,\, \gA)$ classified by $\sU''$, i.e.,~that $2 \cdot \gA(1) = 0$ for every such $(E,\, \gA)$. For this, if $G$ denotes the finite locally free subgroup of order $4$ generated by $\gA$ and $G_2 \subset G$ denotes its standard cyclic subgroup of order $2$, then it suffices to show that the image of $\gA(1)$ in $G/G_2$ vanishes. However, by \cite{KM85}*{6.7.4}, this image generates $G/G_2$, so it remains to recall that, by \cite{KM85}*{6.1.1}, the subscheme of generators 
\[
(G/G_2)^\times \subset G/G_2
\]
is finite locally free of rank $1$ and that, by \cite{KM85}*{12.2.1 and its proof}, $(G/G_2)^\times$ contains the zero section of $G/G_2$ because $G/G_2$ is the Frobenius kernel of $E/G_2$. 
\epf

\brems
\remi \lab{more-U-pr}
We could have analogously constructed a larger $\sU'$ by allowing $\sU''$ to also contain the open locus of $\sY_1(4)_{\bF_2}$ over which the standard cyclic subgroup $\cG_2 \subset \cG$ or order $2$ is of multiplicative type and $\cG/\cG_2$ is \'{e}tale. Over this locus $\cG = \cE[2]$ by \Cref{KM}~\ref{KM-b}, so $2 \cdot \gA(1) = 0$, too. 

\remi
By combining Remarks \ref{comp-mor} and \ref{more-U-pr}, we conclude that the comparison morphism from the coarse moduli space of $\sY_1(4)_{\bF_2}$ towards $Y_1(4)_{\bF_2}$ is finite locally free of rank $2$ over a dense open of the locus over which $\cG$ has a nontrivial multiplicative part and is a universal homeomorphism over the entire $Y_1(4)_{\bF_2}$. In contrast, this morphism is an isomorphism over the locus over which $\cG$ is \'{e}tale because there $\sY_1(4)$ is rigid by \cite{KM85}*{2.7.4 and 1.10.12}.
\erems


\section{Coarse base change fails for $\sX_{\Gamma_1(4) \cap \Gamma_0(4n)}$} \lab{G1-G0}

We seek to explain how a variant of the construction of \S\ref{X1(4)} leads to an infinite family of modular curves that violate coarse base change. We fix any (possibly even) integer $n \in \bZ_{\ge 1}$, set
\[
\Gamma_0(4n) \ce  \{ \p{ \begin{smallmatrix} a & b \\ c & d \end{smallmatrix}  } \in \GL_2(\wh{\bZ})\  |\    c \equiv 0 \bmod 4n  \},
\]
and seek to show that the formation of the coarse moduli space of $\sX_{\Gamma_1(4) \cap \Gamma_0(4n)}$ does not commute with base change to $\bF_2$. The genus of the coarse moduli space of $(\sX_{\Gamma_1(4) \cap \Gamma_0(4n)})_\bQ$ grows unboundedly with $n$ because the same holds already for the coarse moduli space of $(\sX_{\Gamma_0(4n)})_\bQ$.

As in \S\ref{X1(4)}, we first describe $\sY_{\Gamma_1(4) \cap \Gamma_0(4n)}$ in modular terms and explain how we choose $\sU'$.

\bpp[The choice of $\sU'$] \lab{U-prime-15}
By \cite{KM85}*{7.9.4, 7.9.6, and 7.4.2} (with \cite{Ces15a}*{6.3 (a)}),\footnote{We implicitly also use \cite{KM85}*{1.7.2 and 7.3.1} to decompose into primary parts.} the stack $\sY_{\Gamma_1(4) \cap \Gamma_0(4n)}$ parametrizes triples
\[
(E \ra S,\, G,\, \gA\colon \bZ/4\bZ \ra (E/G_n)(S))
\]
consisting of an elliptic curve $E \ra S$ over a variable base scheme $S$, a cyclic $S$-subgroup $G \subset E$ of order $4n$, and a Drinfeld $\bZ/4\bZ$-structure $\gA$ on the $S$-subgroup $G/G_n \subset E/G_n$, where $G_n \subset G$ is the standard cyclic subgroup of order $n$.

We let 
\[
\cE \ra (\sY_{\Gamma_1(4) \cap \Gamma_0(4n)})_{\bF_2}
\]
be the universal elliptic curve in characteristic $2$ and let $\cG \subset \cE$ be the universal $G$. The locus $\sU''$ of $(\sY_{\Gamma_1(4) \cap \Gamma_0(4n)})_{\bF_2}$ over which $\cG$ is of multiplicative type is open, and we let
\[
\sU' \subset \sU'' \qq \text{be the open locus on which the $j$-invariant satisfies} \qq j \neq 0.
\]
\epp

\bthm \lab{fail-XHn}
The $\bF_2$-base change of the coarse space morphism of $\sX_{\Gamma_1(4) \cap \Gamma_0(4n)}$ is not the coarse space morphism of $(\sX_{\Gamma_1(4) \cap \Gamma_0(4n)})_{\bF_2}$. The same holds for any $\bZ_{(2)}$-fiberwise dense open substack of $\sX_{\Gamma_1(4) \cap \Gamma_0(4n)}$.
\ethm

\bpf
Similarly to the proof of \Cref{fail-X1(4)}, it suffices to explain why \Cref{main-crit} applies with 
\[
\sU \ce (\sY_{\Gamma_1(4) \cap \Gamma_0(4n)})_{\bZ_{(2)}} \qq \text{and} \qq \text{$\sU'$ of \S\ref{U-prime-15}.}
\]
The validity of \ref{MC-i} is again supplied by \cite{KM85}*{2.7.4}.

For \ref{MC-ii}, the nonemptiness of the open $\sU''$, and hence also of $\sU'$, follows from \Cref{KM}~\ref{KM-a}, which ensures that the ordinary locus of $\sY(1)_{\bF_2}$ lies in the forgetful image of $\sU''$. Thanks to \Cref{Del}, it remains to show that $2 \cdot \gA(1) = 0$ for every $(E, G, \gA)$ classified by $\sU''$. Since $G/G_{n} \subset E/G_n$ is of multiplicative type and cyclic of order $4$, this follows the same way as in the proof of \Cref{fail-X1(4)}.
\epf


\section{Coarse base change fails for $\sX_1^\bal(4)$} \lab{bal}

In order to complete an open entry in the table \cite{KM85}*{8.5.4}, we investigate coarse base change for one other level. More precisely, in this section we choose the level $H$ to be 
\[
\Gamma_1^\bal(4) \ce \{ \p{ \begin{smallmatrix} a & b \\ c & d \end{smallmatrix}  } \in \GL_2(\wh{\bZ})\  |\  a, d \equiv 1 \bmod 4, \, c \equiv 0 \bmod 4  \}
\]
and seek to show that the formation of the coarse moduli space of 
\[
\sX_1^\bal(4) \ce \sX_{\Gamma_1^\bal(4)}
\] 
does not commute with base change to $\bF_2$. To again apply \Cref{main-crit}, we first explain how we choose $\sU'$.

\bpp[The choice of $\sU'$] \lab{U-prime-2}
By \cite{KM85}*{3.3} (with \cite{KM85}*{7.4.2 (2)} and \cite{Ces15a}*{6.3 (a)}), the stack $\sY_1^\bal(4)$ parametrizes triples
\[
(E \ra S,\, \gA\colon \bZ/4\bZ \ra E(S),\, \gB \colon \bZ/4\bZ \ra (E/G_\gA)(S))
\]
consisting of an elliptic curve $E \ra S$ over a variable base scheme $S$, a Drinfeld $\bZ/4\bZ$-structure $\gA$ on an $S$-subgroup $G_\gA \subset E$, and a Drinfeld $\bZ/4\bZ$-structure $\gB$ on the $S$-subgroup $E[4]/G_\gA \subset E/G_\gA$.

Similarly to \S\ref{U-prime-1}, we consider the universal elliptic curve 
\[
\cE \ra \sY_1^\bal(4)_{\bF_2}
\]
in characteristic $2$ and let $\cG \subset \cE$ be the universal $G_\gA$. The locus $\sU''$ of $\sY_1^\bal(4)_{\bF_2}$ over which the standard cyclic subgroup $\cG_2 \subset \cG$ of order $2$ is of multiplicative type and $\cG/\cG_2$ is \'{e}tale is open, and we let
\[
\sU' \subset \sU'' \qq \text{be the open locus on which the $j$-invariant satisfies} \qq j \neq 0.
\] 
\epp

\bthm \lab{fail-X1bal(4)}
The $\bF_2$-base change of the coarse space morphism of $\sX_1^\bal(4)$ is not the coarse space morphism of $\sX_1^\bal(4)_{\bF_2}$. The same holds for any $\bZ_{(2)}$-fiberwise dense open substack of $\sX_1^\bal(4)$.
\ethm

\bpf
It suffices to explain why \Cref{main-crit} applies with 
\[
\sU \ce \sY_1^\bal(4)_{\bZ_{(2)}} \qq \text{and} \qq \text{$\sU'$ of \S\ref{U-prime-2}.}
\]
The assumption \ref{MC-i} follows from the rigidity of $\sY_1^\bal(4)_{\bZ[\f{1}{2}]}$ supplied by \cite{KM85}*{2.7.4}.

For \ref{MC-ii}, similarly to the proof of \Cref{fail-X1(4)}, we have to show that $\sU'' \neq \emptyset$ and that $\gA(1)$ and $\gB(1)$ are $2$-torsion for every $(E, \gA, \gB)$ classified by $\sU''$. The latter requirement follows from \Cref{KM}~\ref{KM-b}, which ensures that $\cG_{\sU''} = \cE_{\sU''}[2]$. The nonemptiness, on the other hand, follows from \Cref{KM}~\ref{KM-a}, which ensures that the forgetful image of $\sU''$ in $\sY(1)_{\bF_2}$ contains the ordinary locus.
\epf

\end{document}